%% file: main_rev_text_w_supp.tex
\documentclass[aps,pre,showpacs,preprint]{revtex4-2}
\usepackage{graphicx}
\usepackage{amssymb,amsthm,amsfonts,amsmath}
\usepackage{hyperref}

\newcommand{\dd}{\mathrm{d}}

\newcommand{\bbE}{\mathbb{E}}
\newcommand{\bbR}{\mathbb{R}}
\newcommand{\cB}{\mathcal{B}}

\def\Fref#1{Figure~\ref{#1}}
\def\fref#1{Figure~\ref{#1}}

\begin{document}

\title{Quasipotentials for coupled escape problems and the gate-height bifurcation}

\author{Peter Ashwin}
\author{Jennifer Creaser}
\affiliation{Department of Mathematics and Statistics, University of Exeter, Exeter EX4 4QF, UK}
\author{Krasimira Tsaneva-Atanasova}
\affiliation{Department of Mathematics and Statistics, and EPSRC Hub for Quantitative Modelling in Healthcare, University of Exeter, Exeter, EX4 4QJ, UK.}
\affiliation{Institute for Advanced Study, Technical University of Munich, Lichtenbergstrasse 2 a, D-85748 Garching, Germany}

\begin{abstract}
The escape statistics of a gradient dynamical system perturbed by noise can be estimated using properties of the associated potential landscape. More generally, the Freidlin and Wentzell quasipotential (QP) can be used for similar purposes, but computing this is non-trivial and it is only defined relative to some starting point. In this paper we focus on computing quasipotentials for coupled bistable units, numerically solving a Hamilton-Jacobi-Bellman type problem. We analyse noise induced transitions using the QP in cases where there is no potential for the coupled system. Gates (points on the boundary of basin of attraction that have minimal QP relative to that attractor) are used to understand the escape rates from the basin, but these gates can undergo a global change as coupling strength is changed. Such a global gate-height bifurcation is a generic qualitative transitions in the escape properties of parametrised non-gradient dynamical systems for small noise.
\end{abstract}

%\pacs{05.45.Xt (Synchronization; coupled oscillators) 05.40.Ca (Noise)}


\maketitle

\section{Introduction}

Noise induced transitions in dynamical systems are important in understanding a wide variety of phenomena in nonlinear science \cite{kraut2003noise,emenheiser2016patterns,christ2017tristable}, such as in solid state physics \cite{malchow1983stochastic}, brain network dynamics \cite{creaser2018sequential,creaser2020domino} and climate dynamics \cite{lucarini2020global,margazoglou2021dynamical}. For systems that (in the absence of noise) have multiple attractors, the addition of noise will initiate transitions between neighbourhoods of the attractors. These transitions can be understood in terms of ``escape problems'' where one attempts to determine the distribution of times of first exit from a neighbourhood of one attractor $A$ into the basin of another attractor. For low amplitude noise, the exit path from the basin of an attractor $A$ typically goes along a ``most likely path'' (also called instanton) through a ``gate'' on the boundary. The escape rate is then asymptotically given by a large deviation result - there is an ``escape rate'' that scales exponentially with the noise amplitude and the height of the gate.

If the system of interest consists of a network of coupled systems, each of which is multistable, this gives extra structure that one can take advantage of. We consider here a case where individual systems are bistable, but where one of the states (that we call Quiescent) is marginally stable and the other (that we call Active) is substantially more stable, so that the escape rate to return to the Quiescent state is much lower than to escape from the Quiescent state. As noted in \cite{ashwin2017fast,creaser2018sequential}, emergent effects can appear in such sequential escape problems. For a number of uncoupled units, each of which can independently undergo escape, there will clearly be independence of the escape processes, but the presence of coupling can create nontrivial dependence in the escapes and even synchronization of escapes for large enough coupling \cite{berglund2007metastability}. In previous papers \cite{ashwin2017fast,ashwin2018sequential} we highlighted that so-called slow- and fast-domino regimes can arise as the coupling strength changes, in cases where there is escape from a quiescent attractor to an active attractor at rates that are much faster than the reverse escape. The slow-domino regime appears at a critical coupling strength, beyond which escape of one unit to an active state induces the escape of a unit coupled to it, but with some approximately deterministic delay. The fast-domino regime appears at higher couplings and corresponds to cases where the coupling is strong enough such that escape of one unit results in immediate escape of units coupled to it.

These low noise regimes are separated by bifurcations of the basin boundaries of the stable states in the coupled system. Our earlier work \cite{ashwin2017fast} extends the analysis to sequential escapes in a system of three bistable nodes with uni-directional coupling and show how these regimes can be defined and how they affect the escape times and the likely order of escape.

\subsection{The Freidlin-Wentzell quasipotential}

Consider a system on $x\in\bbR^d$ that evolves according to the stochastic differential equation (SDE)
\begin{equation}
    \dd x=f(x)\,\dd t+\sigma \, \dd W_t
    \label{eq:SDE}
\end{equation}
where $f(x)$ is smooth, $W_t$ is a standard Brownian motion on $\bbR^d$ and $\sigma>0$ is a noise amplitude (we assume identity growth in covariance per unit time for simplicity but note that the methodology we use has been extended in \cite{dahiya2018anisotropic} to allow more general matrix-valued $\sigma$). In the limit of $\sigma\rightarrow 0$ one can relate this to the dynamics of the ordinary differential equation (ODE)
\begin{equation}
    \dot x=f(x).
    \label{eq:ODE}
\end{equation}
In the special case, where $f(x)=-\nabla V(x)$ is determined by the gradient of some smooth potential $V:\bbR^d\rightarrow \bbR$, we say (\ref{eq:ODE}) is a gradient system, and in that case one can apply the method of Eyring and Kramers to compute the rate of escape from attractors of (\ref{eq:ODE}) in terms of the potential barrier that must be overcome for a transition to take place (for a review of such methods, see \cite{Berglund2013}).

However, for most choices of $f(x)$, no such potential $V$ exists and so the method above cannot be used - we say the system is non-gradient. Nonetheless, Freidlin and Wentzell \cite{freidlin2012random} introduced a notion of quasipotential (QP) for non-gradient systems, and using this it is possible to obtain rates of escapes and most likely paths taken by escaping trajectories. Similar methodologies has been studied for many years; e.g. \cite{graham1984existence,graham1985weak,zhou2012quasi} but only recently have numerical methods been developed that allow one to explicitly \cite{cameron2012finding} or perturbatively \cite{bouchet2016perturbative} calculate the QP. To define the quasipotential we first need to define the Freidlin-Wentzell action. As explained in \cite{cameron2012finding} the action is a functional depending on a differentiable path $\phi\in C^1([0,T],\bbR^d)$, with $\phi(s)$ in phase space defined for $s\in[0,T]$. It is defined by
\begin{equation}
    S_{T}(\phi)=\frac{1}{2}\int_{0}^{T} \left\|\frac{\dd \phi}{\dd s}(s)-f(\phi(s))\right\|^2\,\dd s.
\end{equation}
Note that $S_{T}(\phi)\geq 0$, and it is zero if and only if $\phi(t)$ is precisely a trajectory $\varphi_t(x_0)$ of (\ref{eq:ODE}) parameterised by time $t\in[0,T]$. We then define the quasipotential with respect to the arbitrary set $A\subset \bbR^d$ as
\begin{equation}
    U_A(x)=\inf\left\{S_{T}(\phi)~:~\phi\in C^1([0,T]),~\phi(0)\in A,~\phi(T)=x,~T>0\right\}.
    \label{eq:UA}
\end{equation}
If $A$ is asymptotically stable then $U_A(x)$ has a minimum at $A$ in a neighbourhood of the basin of $A$.
Note that the quasipotential needs to be defined relative to the subset $A$ of phase space; this is usually chosen to be an attractor. 
Choosing a different subset will give additional information and quasipotentials that differ on different subsets of phase space. 
For any attractor $A$ we write $\cB(A)=\{x\in\bbR^d~:~\varphi_t(x)\rightarrow A\}$ to be the basin of attraction of $A$, where $\varphi_t$ is the flow generated by (\ref{eq:ODE}). 

Computation of the quasipotential is a nontrivial problem in that it requires finding a limiting optimal path $\phi$, and typically there will be non-differentiable points in $U_A(x)$.  The methods that have recently been developed \cite{cameron2012finding,dahiya2018ordered} to compute the quasipotential for low dimensional systems start by transforming to a geometric action, namely \cite[Appendix A]{dahiya2018ordered} shows that
\begin{equation}
    U_A(x)=\inf\left\{ \tilde{S}(\psi)~:~\psi\in C^1([0,L]),~\psi(0)\in A,~\psi(L)=x,~L>0 \right\}.
    \label{eq:UAtilde}
\end{equation}
where we define a geometric action $\tilde{S}$ for the path $\psi(s)$ independent of its parametrization:
\begin{equation}
    \tilde{S}(\psi)=\int_{0}^{L} \|\psi'\| \|f(\psi(s)\|-\psi'.f(\psi(s)) \,\dd s.
\end{equation}
Posing an associated Hamilton-Jacobi-Bellman problem for this geometric action \cite{cameron2012finding} turns the problem of finding $U_A(x)$ from (\ref{eq:UA}) into finding viscosity solutions \footnote{Continuous but not necessarily differentiable solutions} of the following ill-posed Hamilton-Jacobi equation \cite{dahiya2018ordered,yang2019computing}
\begin{equation}
    \|U(x)\|^2+2 f(x)\cdot\nabla U(x)=0,~~U(A)=0.
    \label{eq:HJB}
\end{equation}
This equation is instrumental in finding the minimum action paths that minimise the geometric action relative to some attractor $A$.

Of particular interest is the distribution of first escape times of trajectories $x(t)$ of (\ref{eq:SDE}) from some open set $N$ containing $\cB(A)$ but no other attractors. This is the random variable
$$
\tau_N=\inf \{t>0~:~x(0)\in A~\mbox{and}~x(t)\not \in N\}.
$$
The utility of the quasipotential is that it gives a low-noise asymptotic estimate \cite{freidlin2012random} of the escape time $\tau$ from this neighbourhood $N$ of $\cB(A)$:
\begin{equation}
    \bbE[\tau_{N}] \asymp \exp( U_A(x^*)/\sigma^2)
\end{equation}
as $\sigma\rightarrow 0$, where $x^*$ is a unique point that minimizes $U_A(x)$ for $x\in\partial \cB(A)$. Note that if $x^*$ is a point such that
\begin{equation}
    U_A(x^*)\leq U_A(x),~~x\in \partial \cB(A),
\end{equation}
then we say $x^*$ is the {\em gate} for the basin $\cB(A)$; typically a basin will possess only one gate though this may change as a parameter changes, and there may be multiple gates if there are symmetries of the system that fix the attractor. The relation $\asymp$ indicates logarithmic equivalence; see \cite{Berglund2013,gayrard2004metastability} for precise statements and proofs.

The Ordered Upwind Method (OUM) was introduced in \cite{sethian2001ordered,sethian2003ordered} to approximate solutions of the Hamilton-Jacobi-Bellman equation (\ref{eq:HJB}) on a grid in phase space. This was  subsequently used in \cite{cameron2012finding} to numerically approximate the quasipotential. More recently, this has been improved for 2D phase spaces in \cite{dahiya2018ordered} and we use the latter method. These methods have also been extended to 3D phase spaces in~\cite{yang2019computing,paskal2022efficient} and for anisotropic noise in \cite{dahiya2018anisotropic}. We refer to these papers for more discussion of the algorithms and numerical errors which depend on grid spacing. In our computations we use a $1024\times 1024$ grid of the illustrated part of phase space.

\subsection{Quasipotentials for systems of bistable nodes}

We are not aware of any previous attempts to use quasipotentials to understand cascades of noise-induced escapes for coupled systems. Hence, the aim of this paper is to explore the properties/qualities of a system of coupled nodes using this computational tool. We identify a range of behaviours that are not present in the symmetric/potential case, but do not claim to give an exhaustive theory even in the low noise limit.

\section{Escape for Coupled Bistable Systems}

We consider a network of prototypical bistable nodes governed by the system of SDEs
\begin{equation}
\dd x_i=\bigl[f(x_i,\nu)+\beta\sum_{j\in N_i} (x_j-x_i)\bigr]\dd t + \alpha \,\dd w_i
\label{eq:network}
\end{equation}
where the dynamics of each node is given by
\begin{equation}
\label{eq:onenode}
\dot{x}=f(x,\nu):=-(x-1)(x^2-\nu).
\end{equation}
The coupling strength is $\beta$ and $N_i$ represents the set of neighbours for node $i$.
An independent identically distributed white noise process $\dd w_i$ is added to each node with amplitude $\alpha$. For $0<\nu<1$ the system is bistable with two stable equilibria that we call quiescent ($Q$) and active ($A$) separated by an unstable saddle equilibrium ($S$), we use $\nu=0.01$ and $\alpha = 0.05$ unless otherwise stated. We write these states as $x_Q=-\sqrt{\nu}$, $x_S=\sqrt{\nu}$, $x_A=1$ and note that for small $\nu$ escape from the quiescent state $x_Q$ will be more rapid than from the active state $x_{A}$.

We have used this model with bidirectional coupling to investigate ``domino"-like transitions on small network motifs~\cite{ashwin2017fast}. In this symmetric case the systems can be expressed as a gradient system and the potential landscape $V$ can be computed.

Here we consider the case of two nodes with unidirectional coupling, given by
\begin{align}
dx_1&=\left[f(x_1,\nu)+\beta(x_2-x_1)\right]dt + \alpha\, dw_1,\label{eq:2uni}\\
dx_2&=\left[f(x_2,\nu)\right]dt + \alpha\, dw_2.\nonumber
\end{align}
For chains of nodes, as in \cite{ashwin2017fast} we write $x_{QA}$ to signify states that are continuations from $\beta=0$ of states where $x_1=x_Q$ and $x_2=x_A$ etc. For $\beta>0$ the system is non-gradient and we compute the quasipotential landscape $U$ relative to each attractor $A$. We show how the quasipotential can inform the escape times and escape order for different values of $\beta$. 
We then consider a chain of three nodes, previously considered in \cite{ashwin2017fast}, given by 
\begin{align}
dx_1&=\left[f(x_1,\nu)+\beta(x_2-x_1)\right]dt + \alpha\, dw_1, \nonumber\\
dx_2&=\left[f(x_2,\nu)+\beta(x_3-x_2)\right]dt + \alpha\, dw_2,\label{eq:3uni}\\
dx_3&=\left[f(x_3,\nu)\right]dt + \alpha\, dw_3.\nonumber
\end{align}
We explore what the quasipotential results from two nodes can tell us about cascades of chains of nodes. 

We compare our quasipotential results to numerical simulations of the model computed in {\sc Matlab} using the stochastic Heun method with step size $10^{-3}$.  The initial condition for each realisation is $x_{QQ}$ for the two node system (\ref{eq:2uni}) and $x_{QQQ}$ for the three node system (\ref{eq:3uni}), namely we start with $x_i=x_Q$ for all $i$. We pick a threshold $x_S < \xi\leq x_A$ and compute the time of escape of node $x^{(i)}$ as 
$$
\tau^{(i)}=\inf\{t>0~:~x_i(t)>\xi\}.
$$
We also identify return times, for example the first return to $x_{Q}$ is
\begin{equation}
\tau^{(i)}_R=\inf\{t>0~:~x_i(t)<\xi'~\mbox{ and there is $0<s<t$ with }~x_i(s)>\xi\}
\end{equation}
where $x_Q\leq\xi'< x_S$. Both times $\tau^{(i)}$ and $\tau^{(i)}_R$ are random variables that depend on the coupling strength, the parameters and the particular noise path. Moreover, they only weakly depend on choice of $\xi$ and $\xi'$.
We compute 2000 realisations of the model for each set of parameter values. From this we estimate the mean escape times, mean number of returns and probability of direction of escape.

\begin{figure}[!ht]
\includegraphics[width = 0.5\textwidth]{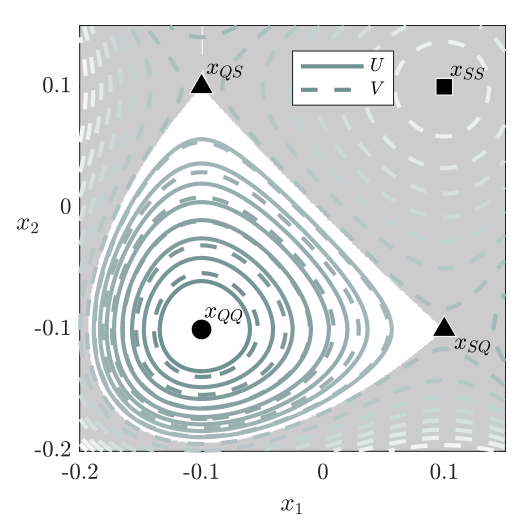} % _2 has mans on
\caption{The potential $V$ (dashed teal) and the scaled quasipotential $U$ (solid teal) computed from $x_{QQ}$ for the uncoupled ($\beta=0$) two node system (\ref{eq:2uni}) with no noise ($\alpha=0$). The white area is the region in which the quasipotential was computed. Equilibria are marked as stable (circle, $x_{QQ}$), saddle (triangle, $x_{SQ}, x_{QS}$) and unstable (square, $x_{SS}$). Note that this symmetric case has two gates: $x_{SQ}$ and $x_{QS}$.}
\label{fig_pot}
\end{figure}

\subsection{Uncoupled}

\begin{figure}[!ht]
\includegraphics[width =\textwidth]{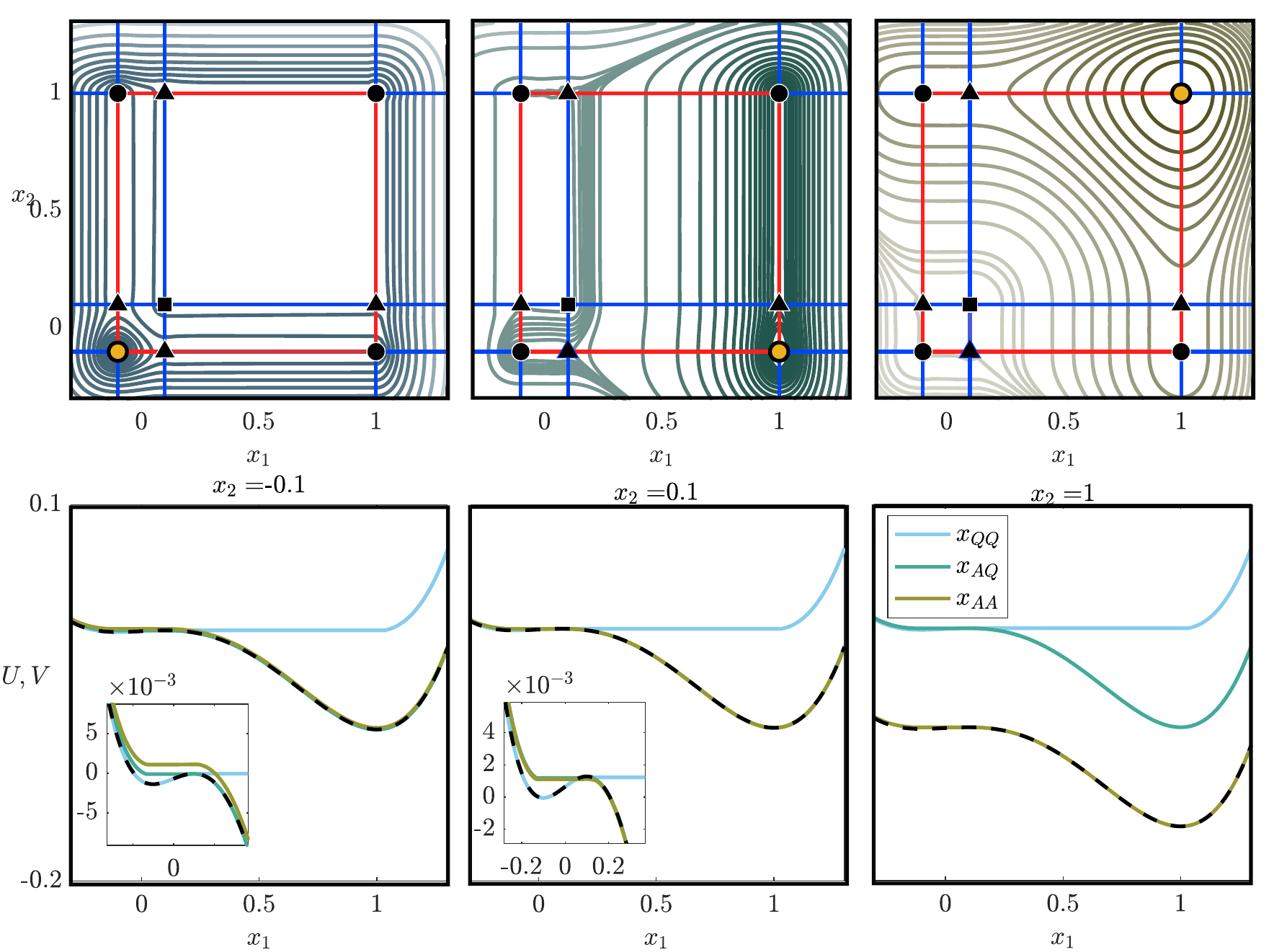}
\caption{Top row: Contour lines of the quasipotential $U$ for (\ref{eq:2uni}) with $\beta=0$ and $\nu=0.01$ computed from $x_{QQ}$ (a), $x_{AQ}$ (b) and $x_{AA}$ (c). Equilibria are marked as in \fref{fig_pot} with the starting point marked in yellow. The red lines show saddle to sink connections, blue lines show other invariant manifolds. The values of the contour lines are chosen for illustrative purposes. Bottom row: The scaled quasipotential plotted against $x_1$ for fixed values of $x_2$ corresponding to the equilibria values, when computed from different starting points. The potential $V$ is shown as a black dashed line, which can be seen to lie on top of one of the coloured lines showing the quasipotentials $U$ from different equilibria. The insets (bottom left and middle) show enlargements of the curves around $x_1=0$; note different values for $U$ corresponding to different attractors.}
\label{fig_qpot}
\end{figure}

We first consider two nodes (\ref{eq:2uni}) with $\beta = 0$ (uncoupled) for which the system admits 9 equilibria corresponding to the states of the system $x_{s_1s_2}$ where $s_i\in\{Q,A,S\}$ is the state of node $i$.
This case is a gradient system for which the potential landscape is
\begin{equation}
V = \frac{x_1^4 + x_2^4}{4} - \frac{x_1^3 + x_2^3}{3} -\nu
        \frac{(x_1^2 + x_2^2)}{2} + \nu(x_1 + x_2).
\label{eq_pot}
\end{equation}
We use the Hamilton-Jacobi-Bellman formulation described above to compute the quasipotential for this system from a given attractor. \Fref{fig_pot} shows the comparison between the contours of the potential landscape $V$ with the (scaled) quasipotential $U$ computed with respect to $x_0 = x_{QQ}$ in a small area of the $(x_1,x_2)$-plane. The quasipotential was computed using the algorithm given in~\cite{dahiya2018ordered} in the white region of the $(x_1,x_2)$-plane. We note that the region in which $U$ is computed appears bounded by the contour line that intersects the two saddle equilibria. In this symmetric case these saddles have the same height in the potential (and quasipotential) landscape. When the level sets of the algorithm reach these gates, it makes an arbitrary choice and continues over one of them. This can be seen in \Fref{fig_qpot} as a thin line along the unstable manifold of the saddle $x_{QS}$ to the computation boundary (the edge of the figure box).
In agreement with the theory, there is a linear relationship $V(x) = U(x)/2 + V(x_0)$ (up to errors from discretization of phase space for computation of $U$), within the basin of attraction of $x_0$ and up to the potential of the gate.

The quasipotential $U$ can be computed starting from any of the four stable equilibria in this system. \fref{fig_qpot} shows the quasipotential computed for the uncoupled system from each of the stable equilibria $x_{QQ}$, $x_{AQ}$ and $x_{AA}$ for the case $\nu=0.01$; due to the symmetry of the uncoupled system we omit $x_{QA}$. The quasipotential can be computed for any arbitrary domain. The contour lines are concentric circles around the starting equilibria up to the nearest saddle, or pair of saddles. From there $U$ deviates from $V$; the quasipotential does not decrease when a saddle or gate is reached, rather it remains constant until the next attractor is reached and only then increases. The most likely path of escape appears as a channel to the next stable point. These channels follow heteroclinic connections from gate to attractor.
The large white regions in the panels correspond to plateaus in the quasipotential. The bottom row of \fref{fig_qpot} shows that this behaviour results in multiple quasipotential height values at each well and gate, depending on which equilibria it is computed from.

\begin{figure}[!ht]
\includegraphics[width =\textwidth]{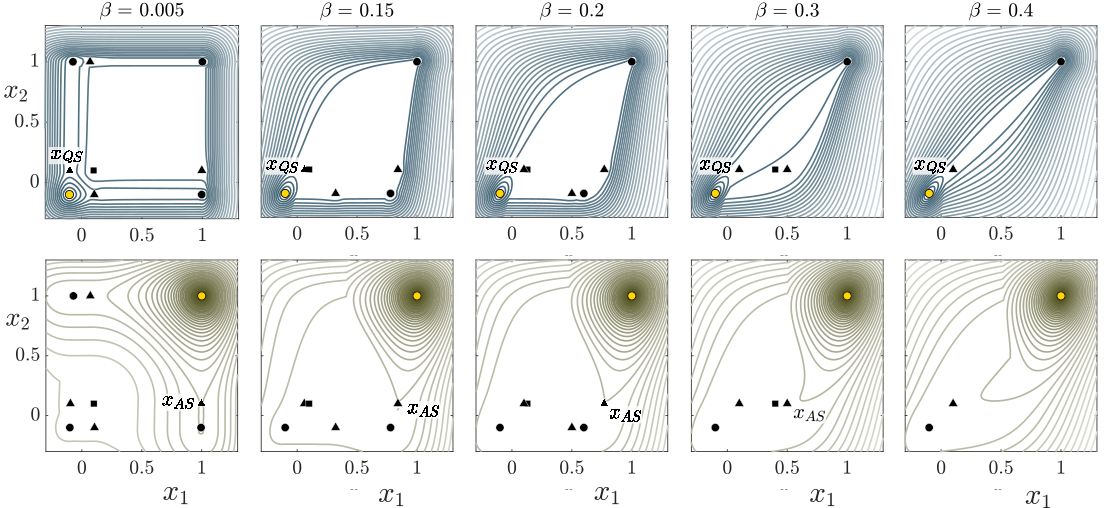}
\caption{Quasipotential landscapes and equilibria for (\ref{eq:2uni}) with $\nu=0.01$ and a range of $\beta$ values. The top row shows the quasipotentials computed from $x_{QQ}$, the bottom row shows quasipotentials computed from $x_{AA}$. Other equilibria are as in \fref{fig_pot}. The gates $x_{QS}$ and $x_{AS}$ are marked to indicate the direction of most likely escape from the starting equilibria (yellow dot).}
\label{fig_qpotuniQQAA}
\end{figure}

\subsection{Uni-directional coupling} 

\begin{figure}[!ht]
\includegraphics[width =\textwidth]{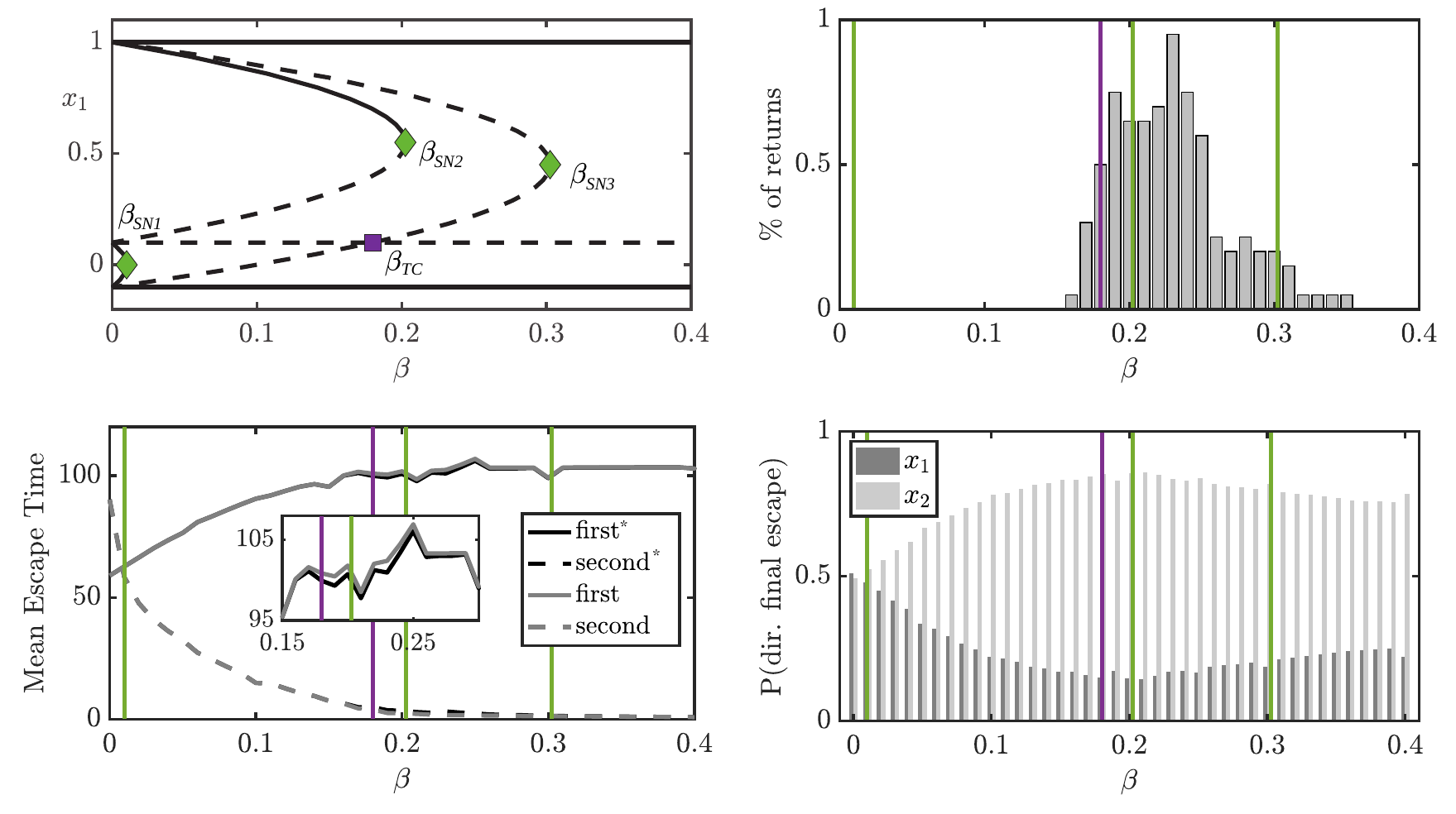}
\caption{Top left: noise-free bifurcation diagram (for $\nu = 0.01$ and $\alpha=0$) for (\ref{eq:2uni}) plotted against $\beta$; green diamonds are saddle-node bifurcations $\beta_{SN\cdot}$, the purple square is the transcritical bifurcation $\beta_{TC}$. The purple and green lines in the remaining panels mark the corresponding bifurcations. 
Top right: the percentage of realisations starting at $x_{QQ}$ that return from $x_{AQ/QA}$ back to $x_{QQ}$ is plotted against $\beta$ for $\alpha=0.05$. 
Bottom left: the mean escape times (detail in inset) for `first*' (the mean first escape from $x_{QQ}$ to $x_{AQ/QA}$), `first' (the mean first escape from $x_{QQ}$ to $x_{AA}$, after possible multiple returns to $x_{QQ}$), `second*' (the mean escape time from $x_{AQ/QA}$ to $x_{AA}$, without passing $x_{QQ}$), `second' (the mean escape time from $x_{AQ/QA}$ to $x_{AA}$, after possible multiple returns to $x_{QQ}$).  
Bottom right: the direction of final escape from the basin of $x_{QQ}$ to $x_{AA}$, plotted against $\beta$. In the uncoupled case ($\beta=0$) the probability is equal at $0.5$. For $\beta>0$ the most likely direction of escape is in the $x_2$ direction. }
\label{fig_times}
\end{figure}

For $\beta>0$ the  system (\ref{eq:2uni}) is non-gradient; we numerically compute the quasipotential $U$ using the method and code presented in \cite{dahiya2018ordered}. \fref{fig_qpotuniQQAA} shows how the equilibria and quasipotential change with the coupling strength $\beta$. The quasipotential is computed from both $x_{QQ}$ and $x_{AA}$.
As $\beta$ increases from $0$ the states $x_{QA}$ and $x_{SA}$ undergo a saddle-node bifurcation of the noise-free system at $\beta_{SN1}=0.01$, denoting the end of the weak coupling regime~\cite{ashwin2017fast}. Unstable states $x_{QS}$ and $x_{SS}$ meet at a transcritical bifurcation at $\beta_{TC} = 0.18$ and states $x_{SQ}$ and $x_{AQ}$ undergo a saddle-node bifurcation at $\beta_{SN2} = 0.2025$. A final saddle-node bifurcation occurs at $\beta_{SN3}=0.3025$ between $x_{SS}$ and $x_{AS}$. The corresponding $\alpha=0$ bifurcation diagram is shown in \fref{fig_times}.

The contours of the quasipotentials shown in \fref{fig_qpotuniQQAA} indicate the global most likely path of escape. From $x_{QQ}$ the most likely path to $x_{AA}$ is via the $x_{QS}$ gate as this gate has the lowest height. This is supported by numerical simulation of the escape times of the two nodes. \fref{fig_times} shows the probability of escape from $x_{QQ}$ is highest in the $x_2$ direction. From $x_{AA}$ the most likely path to $x_{QQ}$ is via the gate $x_{AS}$ for $\beta\leq \beta_{SN3}$ and via $x_{QS}$ for $\beta>\beta_{SN3}$.

\fref{fig_times} also shows that, for certain values of $\beta$, when a realisation escapes in the direction of $x_1$, instead of the more probable $x_2$, there is a chance that it will \emph{return} to $x_{QQ}$ before arriving at $x_{AA}$. This is a deviation from the idealized low noise case where a realization is expected to always follow the most likely path. The proportion of realisations that transition back and forth between (thresholds separating) $x_{QQ}$ and $x_{AQ}$ is around 1\% for the chosen values of $\nu$ and $\alpha$. The mean first and second escape times are also shown in \fref{fig_times} and are $\mathbb{E}[\tau^{(i)}]$.  We distinguish between `first*' escape from the initial condition where no nodes have previously escaped, and `first' where one node has escaped but then the realisation has returned to $x_{QQ}$ before escaping again. Note that this returning behaviour can occur several times before first escape to $x_{AA}$..
The figure also shows times taken for the second node to escape after the first escape, i.e. from $x_{AQ}$ or $x_{QA}$. A low proportion of realisations return, making only a small difference to the mean escape times.

\begin{figure}[!ht]
\includegraphics[width =\textwidth]{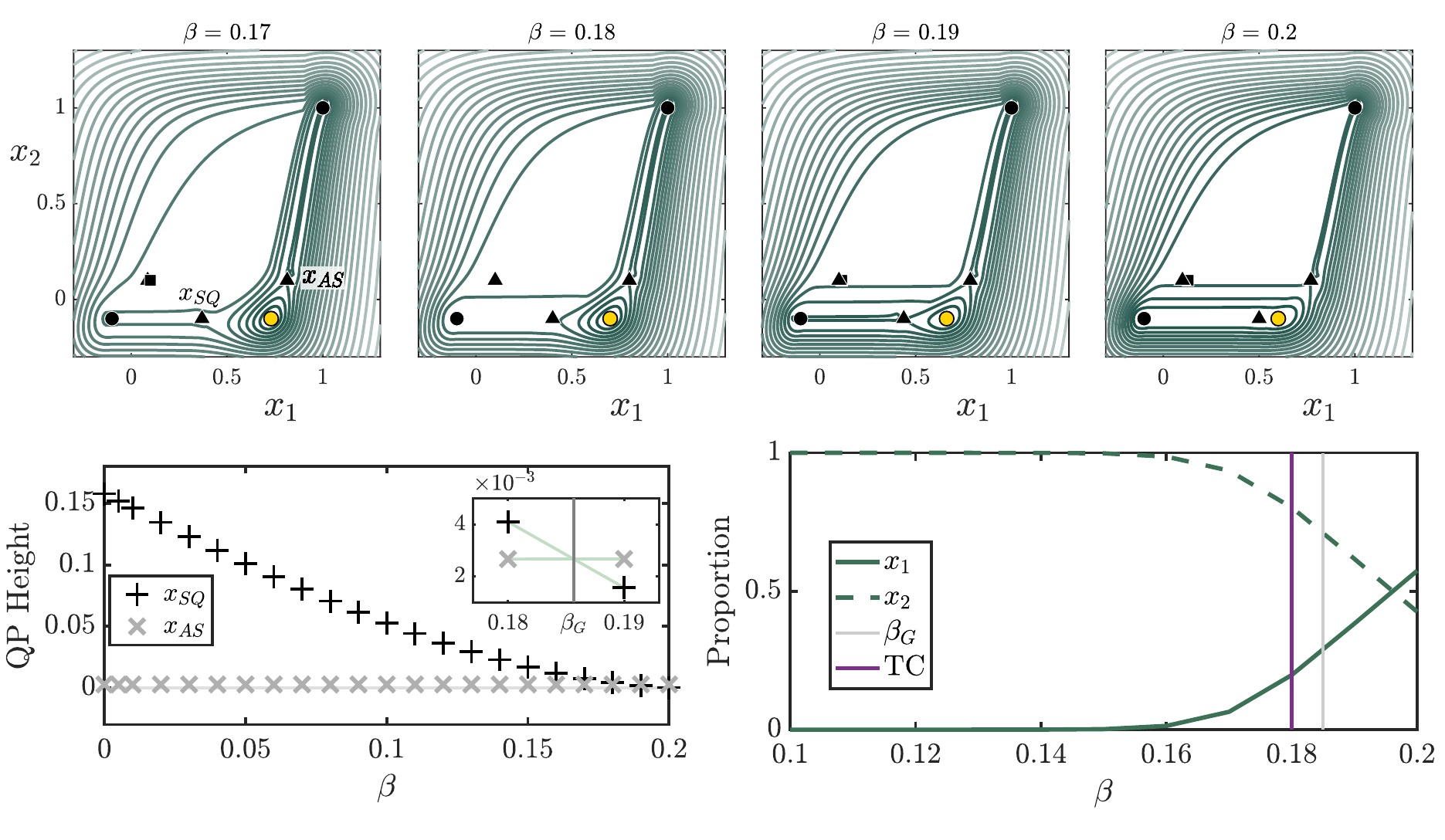}
\caption{The quasipotential for (\ref{eq:2uni}) computed from the equilibria $x_{AQ}$ (yellow dot) for each of the given $\beta$ values (top row); other equilibria are marked as in \fref{fig_pot}.  The height of the quasipotential (QP) at the gates $x_{SQ}$ and $x_{AS}$ against $\beta$ is bottom left. The inset shows the switch in gate heights and the identified \emph{gate-height bifurcation} point $\beta_G \approx 0.1857$. Bottom right shows the proportion of realisations that escape from $x_{AQ}$ in the $x_1$ and $x_2$ directions against $\beta$, the transcritical bifurcation (TC) and the gate-height bifurcation values are marked. The gate-height bifurcation occurs in the $\beta$ range in which there is a local peak in the percentage of returns shown in \fref{fig_times}.}
\label{fig_qpotuniAQ}
\end{figure}

The returning phenomena can be explained by considering the local most likely escape paths from $x_{AQ}$. The quasipotential computed from $x_{AQ}$ is shown for a range of $\beta$ values in \fref{fig_qpotuniAQ}. The top row shows that for $\beta=0.17$ and $\beta=0.18$ there is a channel in the potential landscape between the $x_{AQ}$ state and the full escaped $x_{AA}$ state indicating preferred escape in this direction. Before the second saddle-node at $\beta_{SN2}$ there is a change of preference and the most likely path is back towards $x_{QQ}$. This change of local most likely path can be seen by considering the height of the gates $x_{SQ}$ and $x_{AS}$. The change of likelihood is seen when the height of $x_{SQ}$ becomes lower than $x_{AS}$ for $\beta=\beta_{G}\approx 0.185$. We refer to this qualitative change as a \emph{gate-height bifurcation}. The effect of this bifurcation is observed in the numerical results in \fref{fig_times} at $\beta=0.19$ where there is a local peak in the proportion of returns and a peak in the probability that $x_1$ will be the final direction of escape. Note that the effect of this bifurcation is limited by the saddle-node bifurcation $\beta_{SN2}$ at which $x_{AQ}$ and $x_{SQ}$ coincide.

We note that the quasipotentials computed from $x_{QQ}$ and $x_{AQ}$ look relatively flat between the two equilibria for the $\beta$ values close to the gate-height bifurcation. This allows realisations to have multiple returns, i.e. transition back and forth between the $x_{QQ}$ and $x_{AQ}$ states multiple times before escaping to $x_{AA}$, as observed in the numerical simulations. However, the proportion of realisations with returns is effected by choice of $\alpha$ and $\nu$, as investigated in Appendix~\ref{app:2node}. The proportion of realisations that use the higher gate increases as $\alpha$ increases (for constant $\nu$ and $\beta$), while the proportion of realisations that return decreases as $\nu$ decreases (for constant $\alpha$ and $\beta$).

\section{Quasipotentials for escapes in a chain of bistable systems}

\begin{figure}[!ht]
\includegraphics[width =\textwidth]{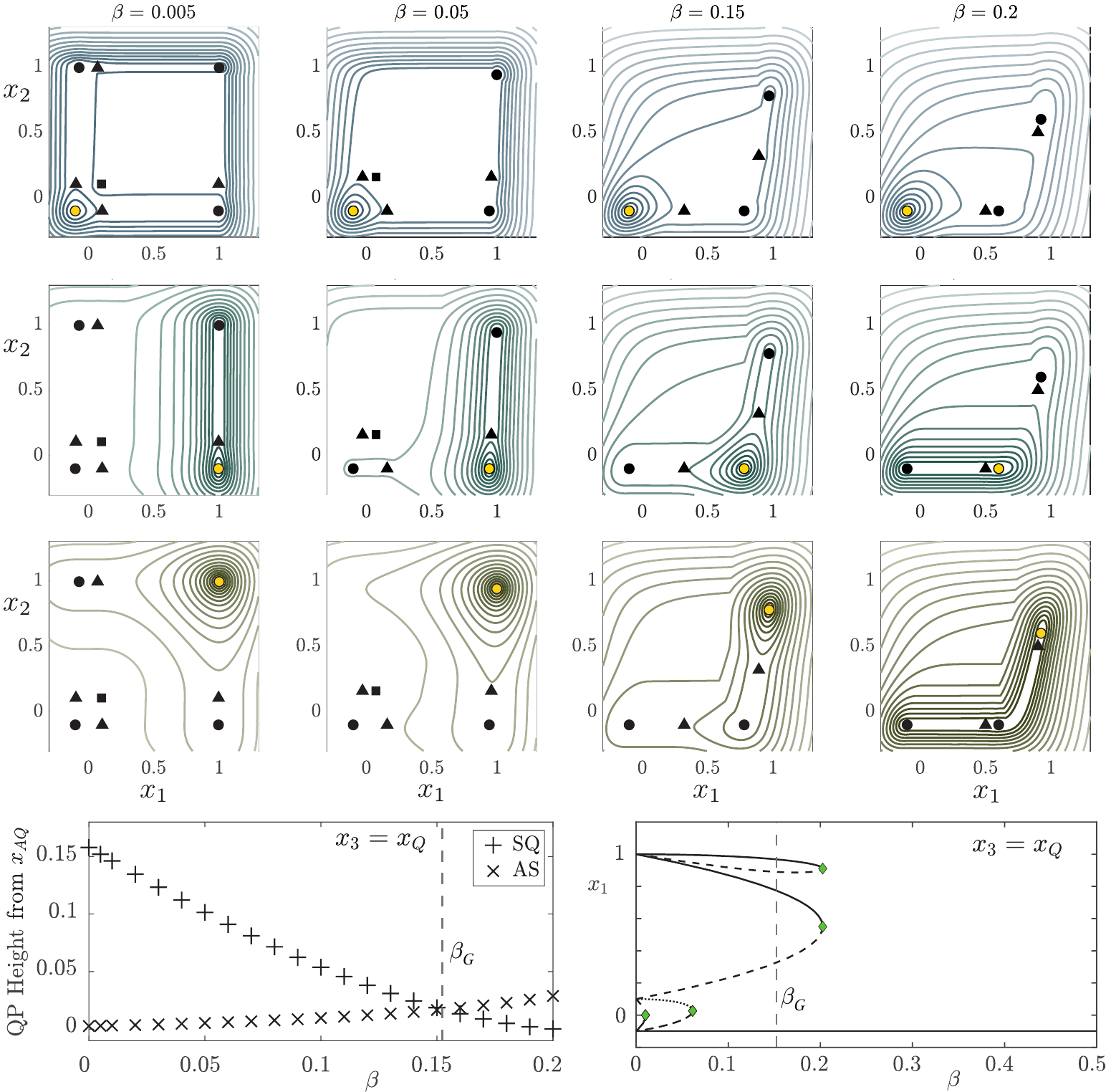}
\caption{Quasipotentials and gate heights for (\ref{eq:3uni}) with $x_3 = x_Q$ for $\nu=0.01$ and varying $\beta$. The yellow dot indicates the starting equilibrium for the quasipotential calculation in each panel; other equilibria are marked as in \fref{fig_pot}. The bifurcation diagram against $\beta$ is shown bottom right with the saddle-node bifurcations marked (green diamonds). The height of the quasipotential computed from $x_{AQ}$ at gates $x_{SQ}$ and $x_{AS}$ is show bottom left; note the exchange of gate heights at $\beta_G\approx 0.1528$.}
\label{fig_qpot3uni}
\end{figure}

We now consider what quasipotentials can tell us about cascades along a chain of three bistable nodes given by system \eqref{eq:3uni}. Although they are not invariant in the presence of noise, one can apply the 2D quasipotential method on planes in the phase space given by $x_3$ fixed at $x_3 = x_Q = -\sqrt{\nu}$ in \eqref{eq:3uni}. This will give an upper bound on the QP in that plane - there may be indirect paths that leave and then return to the plane asymptotically. For $\beta=0$ the equilibria in the system are equivalent to the two node case. The bifurcation diagram against $\beta$, depicted in \fref{fig_qpot3uni}, shows four saddle-node bifurcations. The first at $\beta_{SN1}\approx 0.01$ involves $x_{QAQ}=x_{QA}$ and $x_{SAQ} =x_{SA}$ when compared to the saddles for (\ref{eq:2uni}). The second at $\beta_{SN2}\approx 0.06$ involves unstable states $x_{QSQ}=x_{QS}$ and $x_{SSQ} =x_{SS}$, in contrast to the two node case above. Stable states $x_{AQQ}=x_{AQ}$ and $x_{AAQ} =x_{AA}$ undergo simultaneous saddle-node bifurcations with $x_{SQQ}=x_{SQ}$ and $x_{ASQ} =x_{AS}$, respectively, at $\beta =0.2025$. The only remaining state for $\beta>0.2025$ is $x_{QQQ}=x_{QQ}$.

\fref{fig_qpot3uni} shows the quasipotentials computed from stable states $x_{QQ}$, $x_{AQ}$ and $x_{AA}$ for representative values of $\beta$. The quasipotential computed from $x_{QQ}$ shows the global most likely path to escape is via the $x_2$ direction. From $x_{AA}$ the preferred direction of escape is in the $x_2$ direction towards $x_{AQ}$. To determine preference of direction from $x_{AQ}$ we again consider the height of the gates $x_{SQ}$ and $x_{AS}$. Here the gate-height bifurcation occurs at $\beta=\beta_G \approx 0.1528$ where for $\beta\in[0.15,0.2025]$ the preferred direction is $x_1$ and so to return to $x_{QQ}$. 

\begin{figure}[!ht]
\includegraphics[width =\textwidth]{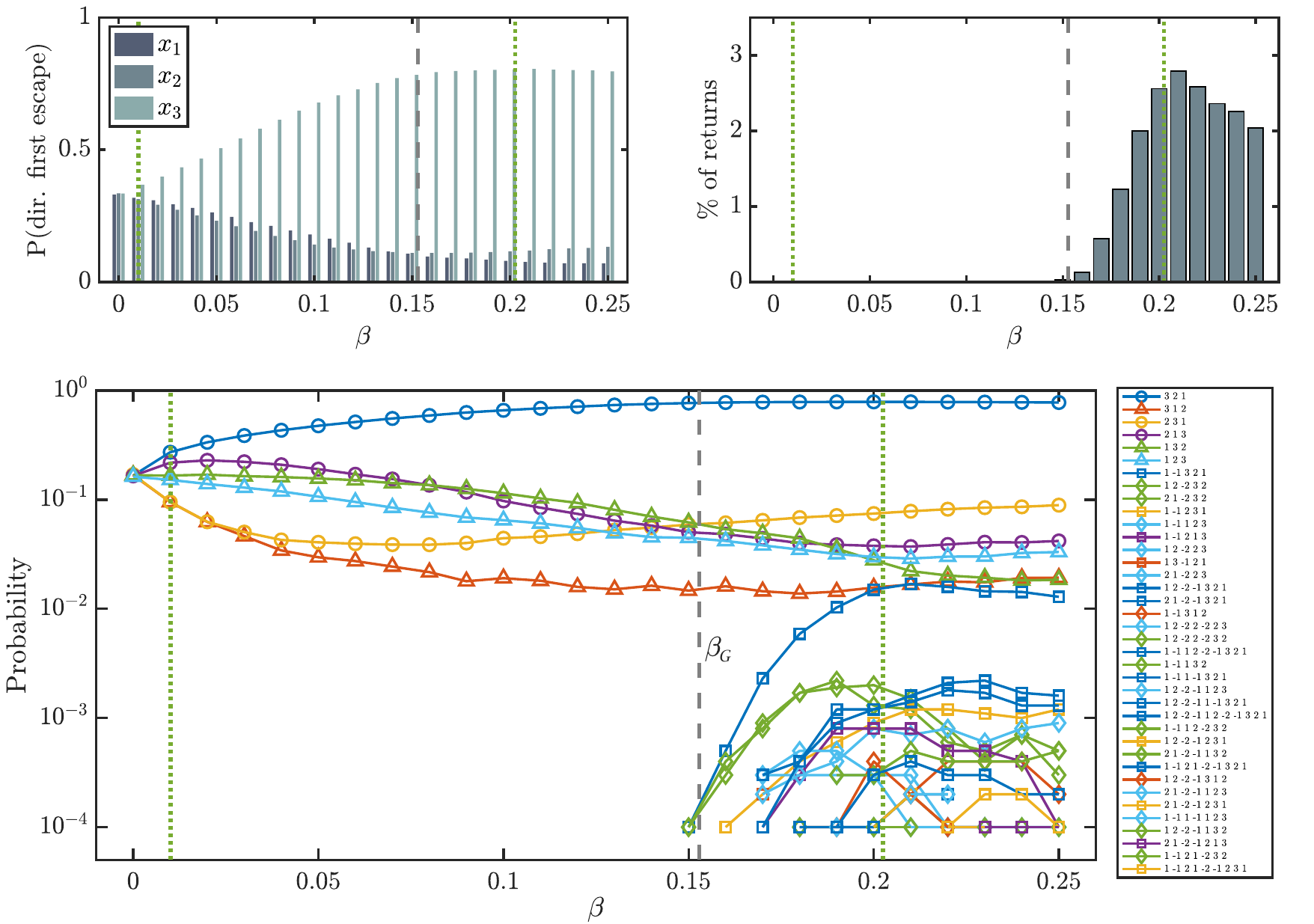}
\caption{Statistics computed from realisations of the full three bistable node system (\ref{eq:3uni}) with $\nu=0.01$, $\alpha=0.05$ and varying $\beta$.  Top row shows the probability of first escape in each direction and the percentage of realisations that return at least once; compare to the two node case \fref{fig_times}. The bottom panel shows the observed sequences and their probabilities where simulations are terminated on the first occasion when all nodes cross a threshold to the state $x_{AAA}$. A returning transition is marked as negative. The sequences with returns are colored according to their final order of escape (cancelling out returns). The circle and square markers denote 2 escapes before 1, and the triangle and diamond markers denote that 1 escapes before 2. Sequences are plotted in order of appearance in increasing $\beta$.
The saddle-node (green dotted lines) and gate-height bifurcations (grey dashed line) are marked in each panel.}
\label{fig_qpot3sim}
\end{figure}

\fref{fig_qpot3sim} shows the direction of escape, percentage of returns and escape sequences computed from the numerical simulations of the full three-node system \eqref{eq:3uni}. Here $10,000$ realisations were computed for $\alpha=0.05$. For $\beta>0$ the most likely direction of first escape is in the $x_3$ direction, and the sequence $[3 2 1]$ is the most likely in line with our previous findings \cite{ashwin2017fast}. The probability of escaping in direction $x_2$ and $x_1$ changes at $\beta_G$ and the percentage of realisations that return is non-zero for $\beta>\beta_G$. The observed sequences are shown with their associated probabilities. The original six sequences (without returns) are all equally probable for $\beta=0$ and no returns are seen. Sequence $[2 3 1]$ initially decreases then increases in probability with increasing $\beta$ and for $\beta>\beta_G$ is the second most likely order. This reflects the change of preference of direction from $x_1$ to $x_2$. This is further supported for large $\beta$ as sequences where 2 escapes before 1 (circle marker, $[2 3 1]$, $[2 1 3]$) become more likely than sequences where 1 escapes before 2 (triangle marker, $[1 3 2]$, $[1 2 3]$).
Sequences with returns appear for $\beta>\beta_G$ and some have multiple returns. They also show that the first and final order of escape differ, for example, for sequence $[1 -1 3 2 1]$ the order of first escape is $[1 3 2]$ but the order of final escape is $[3 2 1]$.  This sequence with one return becomes as probable as some sequences without returns around $\beta = 0.21$. For large $\beta$ returning sequences where the final sequence is $[3 2 1]$ become more probable than other sequences with returns.

For completeness we computed quasipotentials for fixed $x_3=x_{A}$, these are given in Supplementary Material. The stable states $x_{QQA}$, $x_{AQA}$ and $x_{QAA}$ are simultaneously eliminated in saddle-node bifurcations at $\beta=0.01$. There is no gate-height bifurcation and shows that when $x_3$ escapes the other nodes follow almost simultaneously in the so called `fast-domino regime'.

\section{Discussion}

This paper presents the first attempt at computing and analysing quasipotential landscapes of non-gradient systems of coupled bistable nodes. We compute the quasipotentials starting from the meta-stable states and reveal the local and global most likely paths of the system. We identify how these paths change for different values of the coupling strength and through the bifurcations of the noise-free system. 

We introduce the \emph{gate-height bifurcation} of the quasipotential as a global transition for a parameter ($\beta$ in this case) where the values of two local minima $x^a$ and $x^b$ of QP on the basin boundary become equal global minima. In the generic case, the rates of change of the QP with parameter will be unequal meaning that there is transition from one $x^a$ being the gate before the bifurcation to $x^b$ being the gate after the bifurcation. This implies there will be  a qualitative change in the local most likely escape paths from passing through $x^a$ to passing through $x^b$. At the bifurcation there may be more than one likely escape path, each with non-zero probability in the limit $\sigma\rightarrow 0$. Such a gate-height bifurcation allows us to identify regimes where a subset of the realisations of the system are likely to change sequence of visits. In \fref{fig_qpot3sim} we show how this approach gives insight into the timing and order of domino-like cascades of escapes.  We contrast this to the local bifurcation of gates found for two symmetrically coupled bistable units in \cite{ashwin2017fast} and in \cite{creaser2018sequential} where a pitchfork bifurcation of gates distinguishes the slow and fast domino regimes; the degenerate case at bifurcation corresponds to escape over a non-quadratic saddle where a modified version of Kramer's law is needed \cite{berglund2010eyring}.

Our finding that some realisations return to the original state has several implications. It illustrates that the nature of the diffusive coupling depending on the coupling strength can be both activating (or excitatory, i.e. promoting escape to $x_{AA}$) and inhibitory (or suppressing) depending on whether the coupling strength is lower or greater than the gate-height bifurcation value, respectively. An interesting direction for further work would be to investigate how this behaviour depends on the choice of coupling function and whether this affects the robustness of the return of realisations observed with other coupling functions.

From a practical point of view an inhibitory coupling could have implications for preventing or correcting undesirable escape or tipping phenomena. The return of a realisation to its original state indicates that certain  escapes, or tipping events, could be reversed or occur several times before a cascade is triggered.  The standard definition of first escape time and local most likely path should not be considered in isolation in this case. The final escape time, the last escape time of a node given multiple returns, and direction may be more relevant to identify the trigger of the domino effect or cascade. For a system in a given regime, realisations could remain oscillating between two states for a long period. Noisy trajectories have been found to cycle between states in a mean-field model of bursting in neuronal networks~\cite{zonca2022exit}. The authors of \cite{zonca2022exit} compute local potentials for the stable states of this system.  Using quasipotentials to identify the global most likely paths could explain the interplay between escape direction and distributions of escape times in that model.

The coupled bistable model considered here and in \cite{frankowicz1982stochastic} is a simple conceptual example. Its simplicity allows us to compute and analyse the quasipotential landscape for the two and reduced three node examples. The quasipotential approach is widely applicable to analyse transient dynamics in for example, neuroscience~\cite{rabinovich2011robust,zonca2022exit}, gene regulatory networks~\cite{kim2007potential} and climate tipping points~\cite{kronke2020dynamics,margazoglou2021dynamical}. Models for these application areas may include more complex elements in the node dynamics. A natural extension to this work, to make it more applicable to, for example, climate tipping cascades \cite{kronke2020dynamics} would be to consider heterogeneous coupled nodes or, in the case of neuroscience, more physiologically meaningful node dynamics such as those with periodic or excitable dynamics. We also leave for future work investigation of networks of more than three coupled nodes. Note that the quasipotential computation methods \cite{dahiya2018ordered} used here have been extended to stochastic hybrid systems \cite{li2019noise} and to 3D phase spaces in~\cite{yang2019computing,paskal2022efficient}. Explicit computation of the QP in higher-dimensional phase spaces is however challenging - for this reason other methods such as adaptive multilevel splitting \cite{cerou2007adaptive} have been developed to give estimates for large deviation and escape properties in cases where the QP is inaccessible.

\section*{Code availability}
Code for the computations in this paper is available from:

{\tt https://github.com/peterashwin/qp-coupled-escape-2022}.

\section*{Acknowledgements}

We thank Maria Cameron, Valerio Lucarini and Tam\'as T\'el for interesting discussions and advice related to this work. PA and KTA gratefully acknowledge the financial support of the EPSRC via grant EP/T017856/1. PA is partially supported by funding from the European Union's Horizon 2020 research and innovation programme under Grant Agreement 820970 (TiPES) and the UK EPSRC via grant number EP/T018178/1. KTA also acknowledges the support of the Technical University of Munich – Institute for Advanced Study, funded by the German Excellence Initiative.

%------------------------------------------------------------------------------------------------------
%\bibliographystyle{plain}
\bibliographystyle{apsrev4-2}
%\bibliography{qprefs}
%
\input{main_rev_text_w_supp.bbl}

\newpage

\appendix

\section{Two-node system: parameter dependence}
\label{app:2node}

\begin{figure}[!ht]
\includegraphics[width =0.8\textwidth]{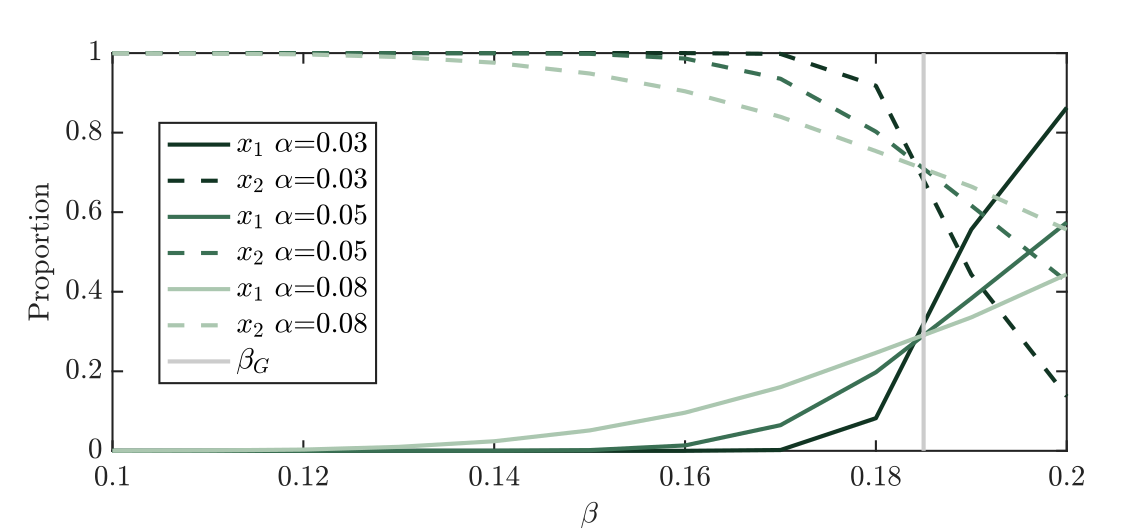}
\caption{Direction of escape from $x_{AQ}$ for $\alpha = {0.03,0.05,0.08}$. The $x_1$ direction is towards $x_{QQ}$ over gate $x_{SQ}$ and the direction $x_2$ is towards $x_{AA}$ over gate $x_{AS}$. The proportion of realisations that escape in the direction of $x_1$ increases as $\alpha$ decreases. with very low noise realisation follow the landscape closely. In the limit at $\alpha\to0$ the proportion of escapes will be 0.5 in each direction at the gate-height bifurcation $\beta_G$.}
\label{fig_alphasEsc}
\end{figure}

\begin{figure}[!ht]
\includegraphics[width =0.9\textwidth]{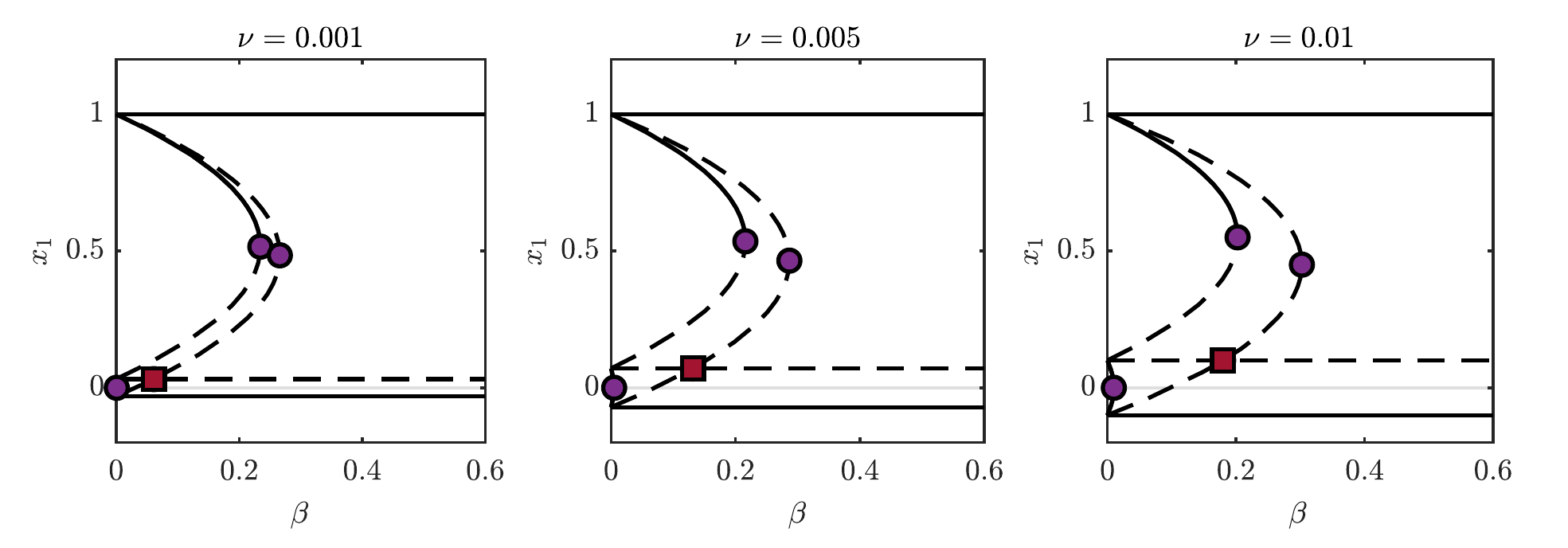}
\includegraphics[width =0.9\textwidth]{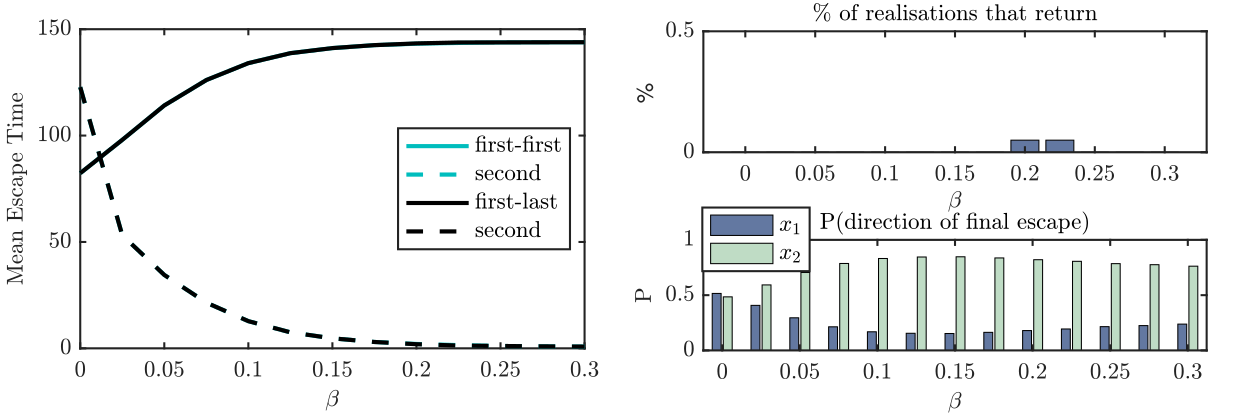}
\caption{Bifurcation diagrams for $\nu = 0.001, 0.005$ and $0.01$. The mean escape times for $\nu=0.005$ are plotted in the bottom left panel, the percentage of realisations with returns and probability of direction of final escape are shown bottom right. Compare to \fref{fig_times}. The percentage of returning realisations is an order of magnitude lower than for $\nu =0.01$ and there is no discernible influence of returns on the mean escape times.}
\label{fig_nus}
\end{figure}

In Figure~\ref{fig_alphasEsc} we illustrate the influence of noise amplitude $\alpha$ on the proportion of escapes to compare with the bottom right panel of Figure~\ref{fig_qpotuniAQ}. 
In \fref{fig_nus} we illustrate the influence of $\nu$ on proportion of escapes to compare with \fref{fig_times}.

%\section{Two-node system: changes with node  parameter}

\clearpage

\section{Three-node system: Node three fixed in active state}
\label{app:3node}

\begin{figure}[!ht]
\includegraphics[width =0.8\textwidth]{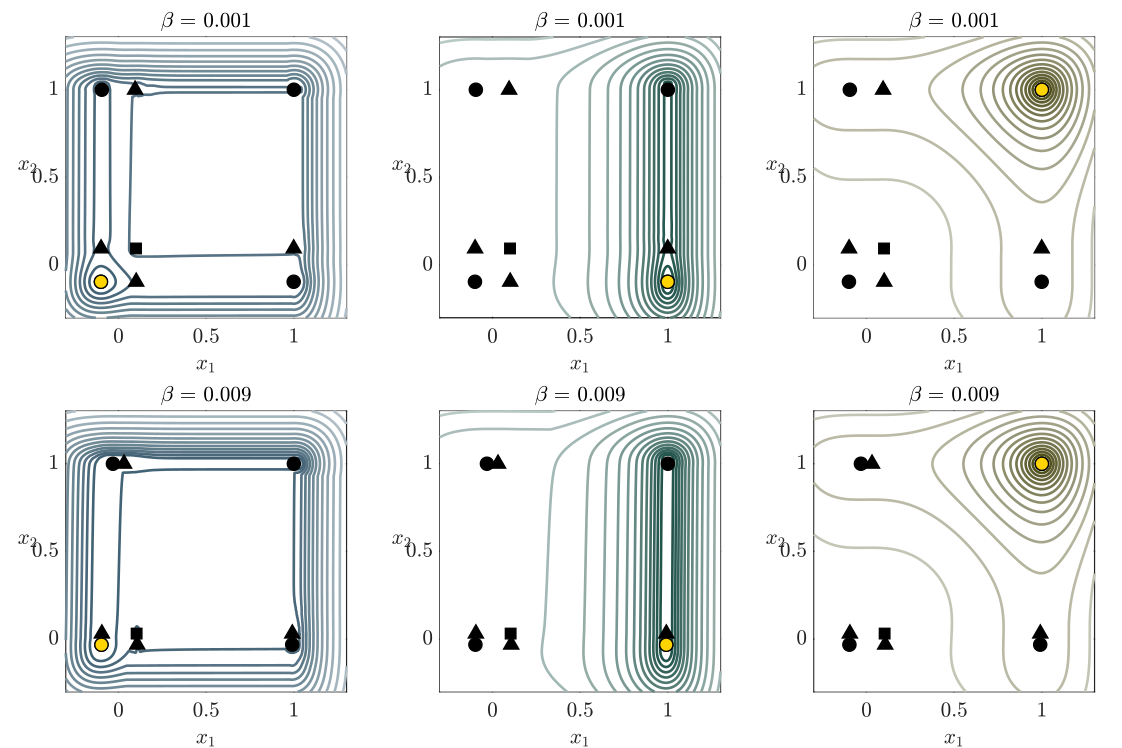}
\caption{Quasipotentials for the three-node system with $x_3 = x_A$ fixed. Potentials are computed from $x_{QQA} = x_{QQ}$ (blue), $x_{AQA} = x_{AQ}$ (green) and $x_{AAA} = x_{AA}$ (dark yellow) the start point is marked as a yellow dot in each panel. The potential from $x_{QAA} = x_{QA}$ is omitted due to the symmetry in the system. Two values of $\beta$ are shown, both are in the weak coupling regime identified in \cite{ashwin2017fast}. For $\beta>0$ the preferred direction of escape from $x_{QQ}$ is $x_2$. From $x_{AQ}$ the preferred direction of escape is  $x_2$, and there is no gate-height bifurcation. For $\beta \approx 0.0101$ there are four simultaneous saddle-node bifurcations of the meta-stable states and for $\beta>0.0101$ only the $x_{AAA}$ state remains.
}
\label{fig_qpot3uniA}
\end{figure}

In \fref{fig_qpot3uniA} we show the QP for the three node system (\ref{eq:3uni}) with $x_3=x_A$ to compare with \fref{fig_qpot3uni}.

\clearpage

\newpage

\section*{Supplementary Material}

\section*{Two-node system: Videos of the two-node potentials}

\begin{itemize}
    \item From QQ: \url{https://youtu.be/JVtxBaydv5Q}
    \item From AA: \url{https://youtu.be/XK6sH8qUHA0}
    \item From AQ: \url{https://youtu.be/00o8Euk9axw}
\end{itemize}

\end{document}